\documentclass[12pt]{article}
\usepackage{amsmath,amsfonts,amssymb,amsthm,amscd}
\title{Remarks on unimodularity}
\date{December 2nd, 2010}
\author{Charlotte Kestner\thanks{Supported by an EPSRC Doctoral Training Grant}\\University of Leeds \and Anand Pillay\thanks{Supported
by EPSRC grant EP/F009712/1}\\University of Leeds}
\newtheorem{Theorem}{Theorem}[section]
\newtheorem{Proposition}[Theorem]{Proposition}
\newtheorem{Definition}[Theorem]{Definition} 
\newtheorem{Remark}[Theorem]{Remark}
\newtheorem{Lemma}[Theorem]{Lemma}

\newtheorem{Problem}[Theorem]{Problem}
\newcommand{\R}{\mathbb R}

\newcommand{\N}{\mathbb N}

\begin{document}
\maketitle

\begin{abstract} 
We clarify and correct some statements and results in the literature concerning unimodularity in the sense of Hrushovski \cite{Hrushovski}, and measurability in the sense of Macpherson and Steinhorn \cite{M-S},  pointing out in particular that the two notions coincide for strongly minimal structures and that another property from \cite{Hrushovski} is strictly weaker, as well as ``completing"  Elwes' proof \cite{Elwes} that measurability implies $1$-basedness for stable theories. 
\end{abstract}

\section{Introduction}
The notion of {\em unimodularity} of a minimal type-definable set $D$ in a stable structure was introduced by Hrushovski \cite{Hrushovski}, in order to to give
a beautiful interpretation and generalization of Zilber's proof that a strongly minimal set $D$ in an $\aleph_{0}$-categorical structure is locally modular. The 
expression {\em unimodularity} refers to the characteristic property that the (locally profinite) group of elementary permutations of any finite-dimensional algebraically closed subset of $D$ is unimodular, namely carries a both left and right invariant Haar measure. On the other hand 
the notion of {\em measurability} of a structure $M$ was introduced by Macpherson and Steinhorn \cite{M-S}, as a direct abstraction of the properties of definable sets in pseudofinite fields obtained by Chatzidakis, van den Dries, and Macintyre \cite{CDM}. The expression {\em measurability} in a mathematical context carries with it a lot of baggage and nuance, so we prefer to refer to this notion of Macpherson and Steinhorn as {\em MS-measurability}. 

Unfortunately what turns out to be a strictly {\em weaker} notion than unimodularity, the non-existence of definable sets $X,Y$, $k\neq \ell$, and definable surjective functions $f:g:X\to Y$ which are $k$-to-$1$, $\ell$-to-$1$, respectively, was claimed in \cite{Hrushovski} to be equivalent to the unimodularity of $D$. We guess this was just an oversight or ``typographical error". But the mistake also surfaces in Elwes' paper \cite{Elwes} as part of a proof that measurable stable structures are $1$-based (see Lemma 6.4 there and its proof), and is repeated in \cite{M-S}  and the survey article \cite{Elwes-Macpherson}. So in the current paper we attempt to clarify the relationships between these various notions and definitions in the stable context, mainly looking at strongly minimal structures. The paper is elementary, and consists mainly of manipulating definitions, and applying results from \cite{Hrushovski} (and \cite{M-S}).

Let us now describe our ``main results" which will appear in section 3. $D$ will denote a strongly minimal structure. Precise definitions will be given in section 2. In Proposition 3.2 we prove the equivalence of  unimodularity and MS-measurability for $D$. This is morally just Lemma 6 of \cite{Hrushovski}, but we also give another equivalent condition involving finite-to-one definable functions. In Proposition 3.2 we give an (easy) example of a weakly unimodular but non unimodular strongly minimal set. In Proposition 3.4 we point out that for  locally modular groups, unimodularity and weak unimodularity coincide.  Finally in Proposition 3.5, we give a correct account of Elwes' theorem from \cite{Elwes} that MS-measurable stable theories are $1$-based. 

Our model-theoretic notation is standard. We assume acquaintance with the basics of geometric stability theory, and the reader is referred to \cite{Pillay-book}, although some of the results we quote are due to Buechler. If (in some ambient model) $a$ is algebraic over $b$, then we write $mult(a/b)$ for the (finite) number of realizations of $tp(a/b)$. Also the expression ``type-definable" refers to definability by a possibly infinite conjunction of formulas over a ``small" set, in an ambient saturated model.

Thanks to Richard Elwes and Dugald Macpherson for several conversations on the topic of the paper, and especially for allowing us to give our commentary on Elwes' results from \cite{Elwes}.

\section{Background on MS-measurability and unimodularity} 
There is a certain tradition in model theory of abstracting features of interesting, concrete examples, to provide a general {\em definition} of a class of structures or theories. What we call here {\em MS-measurability} belongs to this tradition. The definition is simply the conclusion of a theorem from \cite{CDM}. Below we often identify a definable set with the (or a) formula defining it.

\begin{Definition} An $L$-structure $M$ is defined to be MS-measurable if for every nonempty set $X\subseteq M^{n}$ definable (with parameters) in $M$ there is a pair $h(X) = (dim(X),\mu(X))$ with $dim(X)\in \N$ and $\mu(X)\in \R_{>0}$ satisfying the following properties:
\newline
(i) For any $L$-formula  $\phi({\bar x},{\bar y})$,  $\{h(\phi({\bar x},{\bar a})):{\bar a}\in M\}$ is finite and moreover for any given $(k,r)$, 
$\{{\bar a}\in M: h(\phi({\bar x},{\bar a})) = (k,r)\}$ is definable in $M$ without parameters.
\newline
(ii) If $X$ is finite then  $h(X) = (0,|X|)$.
\newline
(iii) If $X,Y \subseteq M^{n}$ are disjoint definable sets, then 
\newline
$dim(X\cup Y) = max\{dim(X),dim(Y)\}$. Also $\mu(X\cup Y)$ equals $\mu(X) + \mu(Y)$ if $dim(X)= dim(Y)$, equals $\mu(X)$ if $dim(X) > dim(Y)$, and equals $\mu(Y)$ if $dim(Y) > dim(X)$.
\newline
(iv)  Suppose $f:X\to Y$ is a definable surjection such that $h(f^{-1}({\bar a}))$ is constant (with value $(d,r)$ say), as ${\bar a}$ varies over $Y$. Then $h(X) = (dim(Y) + d, \mu(Y).r)$. 
\end{Definition}

\begin{Remark} (i) MS-measurability of a structure $M$ is clearly a property of $Th(M)$.
\newline
(ii) If $M$ is MS-measurable by $h = (dim,\mu)$  and $X$ is definable, then $X$ is finite if and only if $dim(X) = 0$.
\end{Remark}

\begin{Lemma} Suppose that $D$ is a strongly minimal structure, which is MS-measurable, witnessed by $h(-) = (dim(-),\mu(-))$. Then for definable sets $X,Y$ in $M$, $dim(X) = dim(Y)$ if and only if $RM(X) = RM(Y)$. In particular MS-measurability of $D$ can also be witnessed by  $(RM(-),\mu(-))$.
\end{Lemma}
\begin{proof} There is no harm in assuming $D$ to be saturated. Suppose that $dim(D) = k$. So $k>0$.  First, from clause (iv) of Definition 2.1 we see that
\newline
{\em Claim.} $dim(D^{n}) = k.n$. 

\vspace{2mm}
\noindent
We now aim to show, by induction on $n$, that for any definable set $X\subseteq D^{n}$, $dim(X) = k.RM(X)$ which suffices. The $n=1$ case is left to the reader. So assume $n > 1$, and $X\subseteq D^{n}$ is definable.
\newline
{\em Case I.}  $RM(X) < n$.
\newline
By compactness we can write $X$ as a disjoint union of of definable sets $X_{1},..,X_{s}$, such that for each $i=1,..,s$ there is $t_{i}\in \N$, a definable subset $Y_{i}$ of $D^{n-1}$ and a $t_{i}$-to-$1$ projection of $X_{i}$ onto $Y_{i}$. We can then apply the induction hypothesis, as well as the definition of MS-measurability, to obtain that $dim(X) = k.RM(X)$ as required.

\vspace{2mm}
\noindent
{\em Case II.}  $RM(X) = n$. 
\newline
Let $Y = D^{n}\setminus X$. So $RM(Y) < n$ and Case I applies. By the Claim and Clause (iii) of Definition 2.1, we conclude that $RM(X) = k.n$. 
\end{proof}

We will make use later of the following result from \cite{M-S}.
\begin{Proposition} For a structure $M$ to be MS-measurable it suffices that there is $h$ satisfying (i)',(ii),(iii), (iv) of Definition 2.1, where (i)' is just (i) for formulas $\phi(x,{\bar y})$ in a single $x$ variable.
\end{Proposition}

We will also use the fact \cite{M-S}  that MS-measurability of $T$ implies MS-measurability of $T^{eq}$. 

\vspace{5mm}
\noindent
We now pass to unimodularity. 
Hrushovski worked with  minimal types in a stable theory, but we  restrict ourselves to strongly minimal sets $D$  (or even strongly minimal structures). So here $D$ denotes a saturated strongly minimal structure.
\begin{Definition}  $D$ is said to be unimodular if for any $n$, if $d_{1},..,d_{n}, d_{1}',..,d_{n}'\in D$, are such that ${\bar d}$ and ${\bar d'}$ are interalgebraic, and $RM(tp({\bar d})) = RM(tp({\bar d'})) = n$ then $mult({\bar d}/{\bar d'}) = mult({\bar d'}/{\bar d})$.
\end{Definition}

The main work in \cite{Hrushovski} was, assuming unimodularity of $D$, to construct a ``Zilber function" on complete types (over parameters), inducing
a Zilber function on type-definable (in particular definable) sets, which, together with Morley rank, turns out to satisfy all properties need for  MS-measurability.  A comprehensive treatment of this also appears in Chapter 2 of \cite{Pillay-book}. The following is the conclusion we need. It appears after Lemma 7 in \cite{Hrushovski}.
\begin{Lemma} Suppose that $D$ is unimodular. Then we can assign, to any type-definable set $X$, a positive real (in fact rational) number $Z(X)$, satisfying the following properties:
\newline
(i) $Z(-)$ is automorphism-invariant,
\newline
(ii) If $X\subseteq Q_{1}\times Q_{2}$ are type-definable sets, and for all $b\in Q_{2}$, $RM(X(b)) = d$ and $Z(X(b)) = r$  (where
$X(b) = \{x\in Q_{1}, (x,b)\in Q\}$), then $Z(Q) = r.Z(Q_{2})$,
\newline
(iii) If $X$ is finite, then $Z(X) = |X|$,
\newline 
(iv) $Z(D) = 1$,
\newline
(v) If $X,Y$ are type-definable disjoint definable subsets of $D^{n}$ then $Z(X\cup Y) = Z(X) + Z(Y)$  if  $RM(X) = RM(Y)$, equals $Z(X)$ if
$RM(X) > RM(Y)$ and equals $Z(Y)$ if $RM(Y) > RM(X)$. 
\end{Lemma}

\vspace{5mm}
\noindent
In the introduction to \cite{Hrushovski}, a strongly minimal structure $D$ was mistakenly defined to be unimodular if whenever $X,Y$ are definable sets in $D$ and $f,g$ are both definable functions from $X$ {\em onto} $Y$, which are $k$-to-$1$, $\ell$ to-$1$ respectively  (with $k,\ell$ positive natural numbers) then $k = \ell$. This property is obviously implied by unimodularity as defined above, using Lemma 2.6 for example. We will call the property {\em weak unimodularity}, although maybe a better expression could be found, and we will see in the next section an easy example of a weakly unimodular but non unimodular strongly minimal set. Note also that this notion of weak unimodularity makes sense in any structure or theory, although again the nomenclature is not ideal as there is not much relation with the locally compact group interpretation.

\section{Results and proofs}

\begin{Proposition} Let $D$ be a strongly minimal set. The following are equivalent:
\newline
(i) $D$ is unimodular,
\newline
(ii) $D$ is MS-measurable,
\newline
(iii) If $X,Y$ are type-definable sets in $D$ and $f,g$ are (type)-definable functions from $X$ onto $Y$ which are $k$-to-$1$, $\ell$-to- $1$, respectively, then  $k = \ell$.
\newline
(iv) Suppose $X,Y$ are definable sets in $D$. Suppose that $U_{1}$, $U_{2}$ are definable subsets of $X$ such that $RM(X\setminus U_{i}) < RM(X)$ for $i=1,2$, and $V_{1}$, $V_{2}$ are definable subsets of $Y$ such that $RM(Y\setminus V_{i}) < RM(Y)$, for $i=1,2$.  Let $f,g$ be definable surjective functions from $U_{1}$ to $V_{1}$, and from $U_{2}$ to $V_{2}$ which are $k$-to-$1$, $\ell$-to-$1$, respectively. Then $k = \ell$.
\end{Proposition}
\begin{proof} As remarked earlier this is essentially the content of Lemma 6 of \cite{Hrushovski}.
\newline 
(i) implies (ii):  Assuming unimodularity of $D$, we will see that defining $h(X) = (RM(X),Z(X))$ for $X$ a definable set in $D$ witnesses 
MS-measurability of $D$. It is well-known that Morley rank in strongly minimal sets satisfies the property: if $f:X \to Y$ is definable and every 
fibre has Morley rank $k$, then $RM(X) = RM(Y) + k$. So using Lemma 2.4, it suffices to check that $h(-) = (RM(-),Z(-))$ satisfies  (i) of Definition 
2.1 for any formula $\phi(x,{\bar y})$ where $x$ is a single variable (ranging over $D$). And this is obvious by strong minimality: there is  $k$ such 
that for any ${\bar b}$, the set defined by $\phi(x,{\bar b})$ is either finite and of cardinality at most $k$ (in which case  $h(\phi(x,{\bar b})) = 
(0,|\phi(x,{\bar b})|)$), or cofinite of co-cardinality at most $k$ (in which case  $h(\phi(x,{\bar b})) = (1,1)$).
\newline
(ii) implies (iv): If $D$ is MS-measurable we saw in Lemma 2.3 that this can be witnessed a by a function $h = (dim,\mu)$ whose ``dimension" component is precisely Morley rank. But then, in the context of (iv) we see that
$\mu(U_{1}) = \mu(U_{2}$ and $\mu(V_{1})= \mu(V_{2})$, whence  (iv) follows immediately.
\newline
(iv) implies (i):  Suppose  ${\bar a}, {\bar b}\in D^{n}$ are interalgebraic, and each of
them realizes $p(x_{1},..,x_{n})$, the ``generic type" of $D^{n}$ (expressing that the $x_{i}$'s are algebraically independent over $\emptyset$). 
Let $k = mult({\bar b}/{\bar a})$ and $\ell = mult({\bar a}/{\bar b})$.  Let 
$q({\bar x},{\bar y}) = tp({\bar a},{\bar b}/\emptyset)$. So $RM(q) = n$ and $q$ has Morley degree $d$ say 
(where possibly $d>1$). We can then find an $L$-formula $\phi({\bar x},{\bar y})\in q$ of Morley rank $n$ and Morley degree $d$, such that $\models \exists^{=k}{\bar y}\phi({\bar a},{\bar y})$ and 
$\models \exists^{=\ell}{\bar x}\phi({\bar x},{\bar b})$. Let $X\subseteq D^{2n}$ be the set defined by $\phi({\bar x},{\bar y})$, $U_{1}$ defined by
$\phi({\bar x},{\bar y})\wedge \exists^{=k}{\bar y}\phi({\bar x},{\bar y})$, and $U_{2}$ defined by 
$\phi({\bar x},{\bar y})\wedge \exists^{=\ell}{\bar x}\phi({\bar x},{\bar y})$. Then clearly $RM(X\setminus U_{i}) < n = RM(X)$ for $i=1,2$.
Let $f:U_{1}\to D^{n}$ be projection to the ${\bar x}$ coordinates, and $g:U_{2}\to D^{n}$ projection to the ${\bar y}$-coordinates. Let $V_{1}\subseteq D^{n}$ be the image of $f$ and $V_{2}\subseteq D^{n}$ the image of $g$. Then each $V_{i}$ has Morley rank $n$ (as it is $\emptyset$-definable and contains a generic point of $D^{n}$), so $RM(D^{n}\setminus V_{i}) < n$. Now $f$ is $k$-to-$1$ and $g$ is $\ell$-to-$1$. So by (iv), $k = \ell$.  We have proved  (i).
\newline
Finally, note that the equivalence of (i) and (iii) is almost immediate. (i) implies (iii) is because of Lemma 2.6. For (iii) implies (i): if 
${\bar d}, {\bar d'}$ are ``generic" interalgebraic tuples from $D^{n}$, then taking $X$ to be the set of realizations of $tp({\bar d},{\bar d'})$ and $Y$ the set of realizations of  $tp({\bar d})$ ($=tp({\bar d'})$), and $f,g$ the obvious projections, then we see that $mult({\bar d}/{\bar d'}) = mult({\bar d'}/{\bar d})$. 
\end{proof}

\begin{Proposition}  There is a weakly unimodular strongly minimal set which is not unimodular.
\end{Proposition}
\begin{proof} This is an easy example. Let the universe of $D$ be $2^{<\omega}$, the collection of finite sequences of $0$'s and $1$'s, equipped with the ``successor relation"  $S$,  where  $S(a,b)$ iff either $b = (a,0)$ or $b = (a,1)$. Then $Th(D)$ has quantifier elimination after adding relations for the compositions $S^{k}$. $Th(D)$ is strongly minimal, and in a saturated model $D'$, if $a$ is ``generic" and $S(a,b)$ holds, then $mult(b/a) =2$ but $mult(a/b) = 1$. A rather painstaking analysis of definable functions will show weak unimodularity of $D$. Note that, working in the standard model $D$, the ``predecessor" function $P$ say, is $2$-to-$1$, its range is all of $D$, but its domain is $D$ with the ``first element" removed. So this predecessor function does not contradict weak unimodularity.
\end{proof}

\begin{Problem}  Is any weakly unimodular  strongly minimal theory  locally modular?

\end{Problem}
\noindent
{\em Comments.} The example in Proposition 3.2 is (geometrically) trivial. On the other hand clearly algebraically closed fields are not weakly unimodular. Conceivably there is a ``Hrushovski construction" of a non locally modular, weakly unimodular strongly minimal set.  However:

\begin{Proposition} Let $D$ be a locally modular strongly minimal group. Then  $Th(D)$ is unimodular iff it is weakly unimodular.
\end{Proposition}
\begin{proof}  We suppose $D$ to be saturated. Assume $D$ is not unimodular.  As  (by Proposition 3.1 for example), unimodularity is independent of adding parameters, we may work over a model, and in particular assume that all types over $\emptyset$ are stationary. So let ${\bar a},{\bar b}$ be tuples in $D^{n}$ witnessing non unimodularity (${\bar a},{\bar b}$ are each generic over $\emptyset$ and $mult({\bar a}/{\bar b}) \neq mult({\bar b}/{\bar a})$). Now by  results in Chapter 4 of \cite{Pillay-book} (and $\omega$-stability), $tp({\bar a},{\bar b})$ is the generic type of a $\emptyset$-definable coset $C$ of a connected $\emptyset$-definable subgroup $H$ of $D^{2n}$. We may assume (by substracting from 
$({\bar a},{\bar b})$ a $\emptyset$-definable point of $C$), that in fact $C = H$. For Morley rank and connectedness reasons, the projections (homomorphisms) $f, g$ of $H$ to the first $n$ coordinates, and the last $n$ coordinates are onto $D^{n}$. Let the ($\emptyset$-definable) kernels of $f,g$ be  $K_{1}, K_{2}$ respectively. As $RM(H) = n$, $K_{1}$ and $K_{2}$ are finite. Let $|K_{1}| = k$ and $|K_{2}| = \ell$. Then $f:H\to D^{n}$ is $k$-to-$1$ and  $g:H\to D^{n}$ is $\ell$-to-$1$. 
\newline
{\em Claim.} $k = mult({\bar b}/{\bar a})$.
\newline
{\em Proof.} If  ${\bar b'}$ realizes  $tp({\bar b}/{\bar a})$  then $({\bar 0},{\bar b'}-{\bar b})\in K_{1}$, so 
$mult({\bar b}/{\bar a}) \leq k$. On the other hand, if $({\bar 0}, {\bar c}) \in K_{1}$, then ${\bar c}\in 
acl(\emptyset)$, so  $RM(tp({\bar a},{\bar b} + {\bar c})) = RM(tp({\bar a},{\bar b})) = n$, whence 
$({\bar a},{\bar b} + {\bar c})$ is also a generic point of $H$ over $\emptyset$, so (by connectedness of $H$) has the same type as $({\bar a},{\bar b})$. Hence $k \leq mult({\bar b}/{\bar a})$, and we get the claim.

\vspace{2mm}
\noindent
Likewise we see that $\ell = mult({\bar a}/{\bar b})$. So by our assumptions, $k\neq \ell$ and $f,g:H\to D^{n}$ witness non weak unimodularity of $D$. 
\end{proof}

\begin{Proposition}  Let $T$ be a (complete) MS-measurable, stable theory. Then $T$ is $1$-based. 
\end{Proposition}
\begin{proof} By Corollary 3.6 of \cite{Elwes-Macpherson}, $T$ is superstable  with finite $R^{\infty}$-rank. By a theorem of Buechler (see Proposition 5.8 of \cite{Pillay-book}), in order to prove that $T$ is $1$-based it suffices to prove that every stationary complete type $p$ of $U$-rank $1$ (in $T^{eq}$) is  locally modular.  By two more results of Buechler (Lemma 3.1 and Proposition 3.2 of \cite{Pillay-book}) we may assume that $p$ has Morley rank $1$ and so is the ``generic type" of some strongly minimal definable set $D$. Now viewing $D$ as a structure in its own right, the measurability of $T^{eq}$ is inherited by $D$ (or $Th(D)$). By Proposition 3.1 and the main theorem of \cite{Hrushovski}, $D$ is locally modular. So $p$ is too. This completes the proof. 

\end{proof}

\end{document}